\documentclass[10pt,psfig]{article}

\usepackage[leqno]{amsmath}
\usepackage{float}
\usepackage{amssymb, amsfonts, amsthm, amscd}
\usepackage{hhline}
\usepackage{array}
\usepackage{pslatex, psfrag, pstricks,epsfig,graphicx}
\usepackage{ifthen}

\newtheoremstyle{dteo}
     {1pt}
     {0pt}
     {\slshape}
     {}
     {\bfseries}
     {:}
     {.5em}
     {}
     
\newtheoremstyle{drem}
     {1pt}
     {1pt}
     {\rmfamily}
     {}
     {\itshape}
     {:}
     {.5em}
     {}

\newtheoremstyle{ddef}
     {1pt}
     {1pt}
     {\rmfamily}
     {}
     {\bfseries}
     {:}
     {.5em}
     {}

\newtheoremstyle{dqst}
     {1pt}
     {1pt}
     {\itshape}
     {}
     {\bfseries}
     {!?!}
     {.5em}
     {}

\theoremstyle{dteo}
\newtheorem{teo}{Theorem}[section]
\newtheorem{lem}[teo]{Lemma}
\newtheorem{cor}[teo]{Corollary}
\newtheorem{prop}[teo]{Proposition}

\theoremstyle{drem}
\newtheorem{rem}[teo]{Remark}
\newtheorem{ex}[teo]{Example}

\theoremstyle{ddef}
\newtheorem{df}[teo]{Definition}

\theoremstyle{dqst}

\def\eqtag{\addtocounter{teo}{1} \tag{\theteo}}

\def\ie{\emph{i.e. }}
\def\eg{\emph{e.g. }}
\def\cf{\emph{cf. }}

\def\eps{\epsilon}
\def\bv{\! \! \mid}
\def\vide{\varnothing}

\def\limm#1{\textrm{\raisebox{.5ex}{\mbox{$\underset{#1}{\lim}$}}} \:}

\def\maxx#1{\textrm{\raisebox{.5ex}{\mbox{$\underset{#1}{\max}$}}} \:}

\def\inff#1{\textrm{\raisebox{.5ex}{\mbox{$\underset{#1}{\inf}$}}} \:}
\def\supp#1{\textrm{\raisebox{.5ex}{\mbox{$\underset{#1}{\sup}$}}} \:}

\def\del{\partial \!}

\def\zz{\mathbb{Z}}
\def\nn{\mathbb{N}}

\def\rr{\mathbb{R}}

\def\dim{\mathrm{dim}}
\def\diam{\mathrm{Diam}}

\def\Id{\mathrm{Id}}

\def\ec{\mathcal{C}}

\def\sm22#1#2#3#4{\left( \begin{smallmatrix}   #1 & #2 \\   #3 & #4 \\  \end{smallmatrix} \right)}

\def\ssi{\Leftrightarrow}
\def\imp{\Rightarrow}

\def\id#1{\mathfrak{#1}}
\def\jo#1{\mathcal{#1}}

\def\nr#1{\left\| #1 \right\|}
\def\pnr#1{\| #1 \|}
\def\abs#1{\left\lvert #1 \right\rvert}

\def\somme#1#2{\overset{#2}{\underset{#1}{\sum}}}

\def\union#1#2{\overset{#2}{\underset{#1}{\cup}}}

\def\wdm{\mathrm{wdim}}
\def\wsp{\mathrm{wspec}}

\def\einj{ \!\! \textrm{ \mbox{ $ {}^{\eps} \!\!\!\!\!\!\! \hookrightarrow$ } } \!\!\! }
\def\frd{\textrm{FilRad}}

\def\moins{\textrm{-}}

\date{November 2, 2007}
\author{Antoine Gournay}
\title{Widths of $l^p$ balls}

\begin{document}

\maketitle

\section{Introduction}

 Let $(X,d)$ be a metric space and  $\eps \in \rr_{>0}$, then we say a map $f:X \to Y$ is an $\eps$-embedding if it is continuous and the diameter of the fibres is less than $\eps$, \ie $\forall y \in Y, \diam f^{-1}(y) \leq \eps$. We will use the notation $f: X \einj Y$. This type of maps, which can be traced at least to the work of Pontryagin (see \cite{Pon} or \cite{HuW}), is related to the notion of Urysohn width (sometimes referred to as Alexandrov width), $a_n(X)$, see \cite{Dra}. It is the smallest real number such that there exists an $\eps$-embedding from $X$ to a $n$-dimensional polyhedron. Surprisingly few estimations of these numbers can be found, and one of the aims of this paper is to present some. However, following \cite{Gro}, we shall introduce:
\begin{df}
$\wdm_\eps X$ is the smallest integer $k$ such that there exists an $\eps$-embedding $f:X \to K$ where $K$ is a k-dimensional polyhedron.
\[
\wdm_\eps (X,d) = \inf_{X \einj K } \dim K.
\]
\end{df}
Thus, it is equivalent to be given all the Urysohn's widths or the whole data of $\wdm_\eps X$ as a function of $\eps$.
\begin{df}
The $\wdm$ spectrum of a metric space $(X,d)$, denoted $\wsp X \subset \zz_{\geq 0} \cup \{+\infty \}$, is the set of values taken by the map $\eps \mapsto \wdm_\eps X$. 
\end{df}
The $a_n(X)$ obviously form an non-increasing sequence, and the points of $\wsp X$ are precisely the integers for which it decreases. We shall be interested in the widths of the following metric spaces: let $B^{l^p(n)}_1$ be the set given by the unit ball in $\rr^n$ for the $l^p$ metric ($\nr{(x_i)}_{l^p} = \big(\sum |x_i|^p\big)^{1/p}$), but look at $ B^{l^p(n)}_1$ with the $l^\infty$ metric (\ie the sup metric of the product). Then
\begin{prop}\label{p1}
$\wsp (B^{l^p(n)}_1, l^\infty ) = \{0,1,\ldots,n\}$, and, $\forall \eps \in \rr_{>0}$,
\[
\wdm_\eps (B_1^{l^p(n)},l^\infty) = \left\{
  \begin{array}{llrcl}
    0  & \textrm{if} &             2 \leq & \eps & \\
    k  & \textrm{if} & 2(k+1)^{-1/p} \leq & \eps & < 2k^{-1/p} \\
    n  & \textrm{if} &                    & \eps &< 2n^{-1/p}
  \end{array}. \right.
\]
\end{prop}
The important outcome of this theorem is that for fixed $\eps$, the $\wdm_\eps (B^{l^p(n)}_1, l^\infty )$ is bounded from below by $\min (n, m(p,\eps) )$ and from above by $\min (n,M(p,\eps))$, where $m,M$ are independent of $n$. As an upshot high values can only be reached for small $\eps$ independantly of $n$. It can be used to show that the mean dimension of the unit ball of $l^p(\Gamma)$, for $\Gamma$ a countable group, with the natural action of $\Gamma$ and the weak-$*$ topology is zero when $p<\infty$ (see \cite{Tsu}). It is one of the possible ways of proving the non-existence of action preserving homeomorphisms between $l^\infty(\Gamma)$ and $l^p(\Gamma)$; a simpler argument would be to notice that with the weak-$*$ topology, $\Gamma$ sends all points of $l^p(\Gamma)$ to $0$ while $l^\infty(\Gamma)$ has many periodic orbits.
\par The behaviour is quite different when balls are looked upon with their natural metric.
\begin{teo} \label{t2}
Let $p \in [1, \infty)$, $n>1$, then $\exists h_n \in \zz $ satisfying $h_n=n/2$ for $n$ even, $h_3=2$ and $h_n = \frac{n+1}{2} $ or $ \frac{n-1}{2}$ otherwise, such that
\[
\{0,h(n),n \} \cup  \subset \wsp (B^{l^p(n)}_1, l^p) \subset \{0 \} \cup (\frac{n}{2}-1,n] \cap \zz.
\]
When $p=2$ or when $p=1$ and there is a Hadamard matrix of rank $n+1$, then $n-1$ also belongs to $\wsp (B^{l^p(n)}_1, l^p)$.
\par More precisely, let $k,n \in \nn$ with $\frac{n}{2}-1<k<n$. Then there exists $b_{n;p} \in [1,2]$  and $c_{k,n;p} \in [1,2)$ such that
\[
\begin{array}{rrlrl}
  \textrm{if} & \eps &\geq 2         & \textrm{then} & \wdm_\eps (B^{l^p(n)}_1, l^p) = 0   \\
  \textrm{if} & \eps &<    2         & \textrm{then} & \wdm_\eps (B^{l^p(n)}_1, l^p) > \frac{n}{2} -1 \\
  \textrm{if} & \eps &\geq c_{k,n;p} & \textrm{then} & \wdm_\eps (B^{l^p(n)}_1, l^p) \leq k  \\
  \textrm{if} & \eps &<    b_{k;p}   & \textrm{then} & \wdm_\eps (B^{l^p(n)}_1, l^p) \geq k  
\end{array}
\]
and, for fixed $n$ and $p$, the sequence $c_{k,n;p}$ is non-increasing. Furthermore, $b_{k;p} \geq 2^{1/p'} \left(1+ \frac{1}{k}\right)^{1/p}$ when $1\leq p \leq 2$, whereas $b_{k;p} \geq 2^{1/p} \left(1+ \frac{1}{k}\right)^{1/p'}$ if $2\leq p < \infty$. 
\par Additionally, in the Euclidean case ($p=2$), we have that $b_{n;2} = c_{n-1,n;2} = \sqrt{2 (1+\frac{1}{n})}$, while in the $2$-dimensional case $b_{2;p} \geq \max (2^{1/p},2^{1/p'})$ for any $p \in [1,\infty]$. Also, if $p=1$, and there is a Hadamard matrix in dimension $n+1$, then $b_{n;1} = c_{n-1,n;1} = \left(1+ \frac{1}{n}\right)$. Finally, when $n=3$, $\forall \eps>0, \wdm_\eps B^{l^p(n)}_1 \neq 1$ and $c_{2,3;p} \leq 2(\frac{2}{3})^{1/p}$, which means in particular that $c_{2,3;p} = b_{3;p}$ when $p\in [1,2]$.
\end{teo}
Various techniques are involved to achieve this result; they will be exposed in section \ref{spt2}. While upper bounds on $\wdm_\eps X$ are obtained by writing down explicit maps to a space of the proper dimension (these constructions use Hadamard matrices), lower bounds are found as consequences of the Borsuk-Ulam theorem, the filling radius of spheres, and lower bounds for the diameter of sets of $n+1$ points not contained in an open hemisphere (obtained by methods very close to those of \cite{PI}). We are also able to give a complete description in dimension $3$ for $1 \leq p \leq 2$. 

\section{Properties of $\wdm_\eps$}
Here are a few well established results; they can be found in \cite{Coo}, \cite{CK}, \cite{Lin}, and \cite{LW}.
\begin{prop} \label{propbase} \renewcommand{\labelenumi}{{\normalfont \alph{enumi}.}}
Let $(X,d)$ and $(X',d')$ be two metric spaces. $\wdm_\eps$ has the following properties:
  \begin{enumerate}
    \item If $X$ admits a triangulation, $\wdm_\eps (X,d) \leq \dim X$.
    \item The function $\eps \mapsto \wdm_\eps X $ is non-increasing.
    \item Let $X_i$ be the connected components of $X$, then $\wdm_\eps (X,d) =0 \ssi \eps \geq \maxx{i} \diam X_i$. 
    \item If $f:(X,d) \to (X',d')$ is a continuous function such that $d(x_1,x_2) \leq C d'(f(x_1),f(x_2))$ where $C \in ]0,\infty[$, then $\wdm_\eps (X,d) \leq \wdm_{\eps/C} (X',d')$.
    \item Dilations behave as expected, \ie let $f:(X,d) \to (X',d') $ be an homeomorphism such that $d(x_1,x_2) = C d'(f(x_1),f(x_2))$; this equality passes through to the $\wdm$: $\wdm_\eps (X,d) = \wdm_{\eps/C} (X',d') $.
    \item If $X$ is compact, then $\forall \eps >0, \wdm_\eps (X,d) < \infty$.
  \end{enumerate} 
\end{prop}
\begin{proof} \renewcommand{\labelenumi}{{\normalfont \alph{enumi}.}}
They are brought forth by the following remarks:
  \begin{enumerate}
  \item If $\dim X =\infty$, the statement is trivial. For $X$ a finite-dimensional space, it suffices to look at the identity map from $X $ to its triangulation $ T(X)$, which is continuous and injective, thus an $\eps$-embedding $\forall \eps $.
  \item If $\eps \leq \eps'$, an $\eps$-embedding is also an $\eps'$-embedding.
  \item If $\wdm_\eps X =0 $ then $\exists \phi:X \einj K$ where $K$ is a totally discontinuous space. $\forall k \in K, \phi^{-1}(k)$ is both open and closed, which implies that it contains at least one connected component, consequently $\diam X_i \leq \eps$. On the other hand, if $\eps \geq \diam X_i$ the map that sends every $X_i$ to a point is an $\eps$-embedding.
  \item If $\wdm_{\eps/C} X' = n$, there exists an $\frac{\eps}{C}$-embedding $\phi:X' \to K$ with $\dim K=n$. Noticing that the map $\phi \circ f$ is an $\eps$-embedding from $X$ to $K$ allows us to sustain the claimed inequality.
  \item This statement is a simple application of the preceding for $f$ and $f^{-1}$.
  \item To show that $\wdm_\eps$ is finite, we will use the nerve of a covering; see \cite[\S V.9]{HuW} for example. Given a covering of $X$ by balls of radius less than $\eps/2$, there exists, by compactness, a finite subcovering. Thus, sending $X$ to the nerve of this finite covering is an $\eps$-immersion in a finite dimensional polyhedron.\hfill \qedhere
  \end{enumerate}
\end{proof}
\indent Another property worth noticing is that $\limm{\eps \to 0} \wdm_\eps (X,d) = \dim X$ for compact $X$; we refer the reader to \cite[prop 4.5.1]{Coo}. Reading \cite[app.1]{Gfrm} leads to believe that there is a strong relation between $\wdm$ and the quantities defined therein ($\textrm{Rad}_k$ and $\diam_k$); the existence of a relation between $\wdm$ and the filling radius becomes a natural idea, implicit in \cite[1.1B]{Gro}. We shall make a small parenthesis to remind the reader of the definition of this concept, it is advised to look in \cite[\S 1]{Gfrm} for a detailed discussion. 
\par Let $(X,d)$ be a compact metric space of dimension $n$, and let $L^\infty (X)$ be the (Banach) space of real-valued bounded functions on $X$, with the norm $\pnr{f}_{L^\infty}= \supp{x\in X} \abs{f(x)}$. The metric on $X$ yields an isometric embedding of $X$ in $L^\infty (X)$, known as the Kuratowski embedding:
\[
\begin{array}{rcl}
I_X: X &\to    & L^\infty(X) \\
   x   &\mapsto& f_x(x') = d(x,x').
\end{array}
\]
The triangle inequality ensures that this is an isometry:
\[
\nr{f_x - f_{x'}}_{L^\infty} = \supp{x'' \in X} \abs{d(x,x'')-d(x',x'')} = d(x,x').
\]
Denote by $U_\eps(X)$ the neighborhood of $X \subset L^\infty(X)$ given by all points at distance less than $\eps$ from $X$,
\[
\ie \qquad U_\eps(X) = \big\{f \in L^\infty(X) \big| \inff{x \in X} \pnr{f-f_x}_{L^\infty} < \eps  \big\} .
\]
\begin{df} 
The filling radius of a $n$-dimensional compact metric space $X$, written $\frd X$, is defined as the smallest $\eps$ such that $X$ bounds in $U_\eps(X)$, \ie $I_X(X) \subset U_\eps(X)$ induces a trivial homomorphism in simplicial homology $H_n(X) \to H_n(U_\eps(X))$.
\end{df}
Though $\frd$ can be defined for an arbitrary embedding, we will only be concerned with the Kuratowski embedding.
\begin{lem} \label{kitu}
Let $(X,d)$ be a $n$-dimensional compact metric space, $k<n$ an integer, and $Y\subset X$ a $k$-dimensional closed set representing a trivial (simplicial) homology class in $H_k(X)$. Then
\[ \eps < 2 \frd Y \imp \wdm_\eps (X,d) > k. \]
If we remove the assumption that $[Y] \in H_k(X)$ be trivial, the inequality is no longer strict: $\wdm_\eps (X,d) \geq k$.
\end{lem}
\begin{proof}
Let us show that $\wdm_\eps (X,d) \leq k \imp \eps \geq 2 \frd Y$. Given an $\eps$-embedding $\phi: X \einj K$, then $\phi(Y) \subset K$ bounds, since $\phi_*[Y]=0$ as $[Y]=0$ in $H_k(X)$. Since $\dim K \leq k = \dim Y$, the chain representing $\phi(Y)$ is trivial. Compactness of $X$ allows us to suppose that $\phi$ is onto a compact $K$. Otherwise, we restrict the target to $\phi(X)$. We will now produce a map $Y \to L^\infty(Y)$ whose image is contained in $U_{\eps/2}(Y)$, so that $Y$ will bound in its $\frac{\eps}{2}$-neighborhood. This will mean that $\eps \geq 2\frd Y$. Let
\[
\begin{array}{ccc}
\begin{array}{rcl}
Q: K &\to    & L^\infty(X) \\
   k &\mapsto& g_k(x'') = \eps/2 + \inff{x' \in \phi^{-1}(k)} d(x'',x')
\end{array}
& \textrm{, and}&
\begin{array}{rcl}
\rho_Y: L^\infty(X) &\to    & L^\infty(Y) \\
       f          &\mapsto& f\bv_Y.
\end{array}
\end{array}
\]
First, notice that $\rho_Y \circ Q \circ \phi (Y) \subset U_{\delta+\eps/2}(Y), \forall \delta>0$ :
\[
\nr{\rho_Y \circ Q \circ \phi (y) -I_Y(y) }_{L^\infty} 
   = \supp{y'' \in Y} \abs{\frac{\eps}{2} + \bigg[ \inff{y' \in \phi^{-1}(\phi(y))} d(y'',y') \bigg] - d(y'',y)} 
   = \eps/2,
\]
since $\phi$ is an $\eps$-embedding. Second, $(\rho_Y \circ Q \circ \phi)_* [Y]=0 $ and $(\rho_Y \circ Q \circ \phi) \sim I_Y$ in $U_{\delta+\eps/2}(Y)$, as $L^\infty(Y)$ is a vector space. Consequently, $[I_Y(Y)]=0$ and $\eps \geq 2 \frd Y$, by letting $\delta \to 0 $. 
\par If $[Y] \neq 0 \subset H_k(X)$, the proof still follows by taking $K$ of dimension $k-1$: the homology class $\phi_*[Y]$ is then inevitably trivial, since $K$ has no rank $k$ homology.
\end{proof}
\par Thus, calculating $\frd$ is a good starting point. The following lemma consists of a lower bound for $\frd$: 
\begin{lem} \label{wdmineq}
Let $X$ be a closed convex set in a $n$-dimensional normed vector space. Suppose it contains a point $x_0$ such that the convex hull of $n+1$ points on $\del X$ whose diameter is $<a$ excludes $x_0$. Then $\frd \del X \geq a/2 $, and, using lemma \ref{kitu}, $\eps < a \imp \wdm_\eps X =n$.
\end{lem}
\begin{proof}
Suppose that $Y=\del X$ has a filling radius less than $a/2$. Then, $\exists \eps>0$ and $\exists P$ a polyhedron such that $Y$ bounds in $P$, $P\subset U_{\frac{a}{2}-\eps}(Y)$ and that the simplices of $P$ have a diameter less than $\eps$. To any vertex $p \in P$ it is possible to associate $f(p) \in I_Y(Y)$ so that $\nr{p,f(p)}_{L^\infty(Y)}<\frac{a}{2}-\eps$. Let $p_0, \ldots, p_n$ be a $n$-simplex of $P$, 
\[
\diam \{ f(p_0), \ldots, f(p_n) \} < 2(\frac{a}{2}-\eps)+\eps <a-\eps<a.
\]
Since $I_Y$ is an isometry, $f(p_i)$ can be seen as points of $Y$ without changing the diameter of the set they form. The convex hull of these $f(p_i)$ in $B$ will not contain $x_0$: their distance to $f(p_0)$ is $<a$ which excludes $x_0$. Let $\pi$ be the projection away from $x_0$, that is associate to $x\in X$, the point $\pi(x) \in \del X$ on the half-line joining $x_0$ to $x$. Using $\pi$, the $n$-simplex generated by the $f(p_i)$ yields a simplex in $Y$. Thus we extended $f$ to a retraction $r$ from $P$ to $Y$. Let $c$ be a $n$-chain of $P$ which bounds $Y$, \ie $[Y]=\delta c$. A contradiction becomes apparent: $[Y]=r_*[Y]=r_* \delta c = \delta r_* c$. Indeed, if that was to be true, $Y$, which is $n-1$ dimensional would be bounding an $n$-dimensional chain in $Y$. Hence $\frd Y>a/2$.
\end{proof}
This yields, for example:
\begin{lem} \label{lemgro} 
(\cf \cite[1.1B]{Gro}) Let $B$ be the unit ball of a $n$-dimensional Banach space, then $\forall \eps<1, \wdm_\eps B = n$.
\end{lem}
\begin{proof}
Any set of $n+1$ points on $Y=\del B$ whose diameter is less than $1$ does not contain the origin in its convex hull. So according to lemma \ref{wdmineq}, $\frd Y>1/2$, and since $Y$ is a closed set of dimension $n-1$ whose homology class is trivial in $B$, we conclude by applying lemma \ref{kitu}.
\end{proof}
Let us emphasise this important fact on $l^\infty$ balls in finite dimensional space.
\begin{lem} \label{wdmli}
Let $B^{l^\infty(n)}_1 = [-1,1]^n$ be the unit cube of $\rr^n$ with the product (supremum) metric, then
\[
\wdm_\eps B^{l^\infty(n)}_1 = \left\{
  \begin{array}{llc}
    0  & \textrm{if} &                       \eps \geq 2 \\
    n  & \textrm{if} &                       \eps < 2
  \end{array}. \right.
\]
\end{lem}
This lemma will be used in the proof of proposition \ref{p1}. Its proof, which uses the Brouwer fixed point theorem and the Lebesgue lemma, can be found in \cite[lem 3.2]{LW}, \cite[prop 2.7]{CK} or  \cite[prop 4.5.4]{Coo}.
\begin{proof}[Proof of proposition \ref{p1}: ]
We first show the lower bound on $\wdm_\eps$. In a $k$-dimensional space, the $l^\infty $ ball of radius $k^{-1/p}$ is included in the $l^p$ ball: $B^{l^\infty(k)}_{k^{-1/p}} \subset B^{l^p(k)}_1$, as $\nr{x}_{l^p(k)} \leq k^{1/p} \nr{x}_{l^\infty(k)}$. Since $B^{l^p(k)}_1 \subset B^{l^p(n)}_1$, by \ref{propbase}.d, we are assured that, if $ B^{l^p(n)}_1$ is considered with the $l^\infty$ metric, $\eps < 2k^{-1/p}$ implies that $\wdm_\eps (B^{l^p(n)}_1,l^\infty) \geq k$.
\par To get the upper bound, we give explicit $\eps$-embeddings to finite dimensional polyhedra. This will be done by projecting onto the union of $(n-j)$-dimensional coordinates hyperplanes (whose points have at least $j$ coordinates equal to $0$). Project a point $x \in B^{l^p(n)}_1$ by the map $\pi_j$ as follows: let $m$ be its $j^{th}$ smallest coordinate (in absolute value), set it and all the smaller coordinates to $0$, other coordinates are substracted $m$ if they are positive or added $m$ if they are negative.
\par Denote by $\vec{\eps}$ an element of $\{-1,1\}^n$ and $\vec{\eps}_{\setminus A}$ the same vector in which $\forall i \in A$, $\eps_i$ is replaced by $0$. The largest fibre of the map $\pi_j$ is 
\[ 
\pi_j^{-1}(0) = \union{\vec{\eps}, i_1, \ldots ,i_{j-1}}{} \{ \lambda_0 \vec{\eps} + \somme{1\leq l\leq j-1}{} \lambda_l \vec{\eps}_{\setminus \{i_1, \ldots, i_l \} } | \lambda_i \in \rr_{\geq 0}\} \cap B^{l^p(n)}_1.
\]
Its diameter is achieved by $s_0 = \big( (n-j+1)^{-1/p}, \ldots , (n-j+1)^{-1/p},0, \ldots,0 \big)$ and $-s_0$; thus $\diam \pi_j^{-1}(0) = 2(n-j+1)^{-1/p}$. $\pi_j$ allows us to assert that $\eps >2(n-j+1)^{-1/p} \imp \wdm_\eps (B^{l^p(n)}_1,l^\infty) \leq n-j $, by realising a continuous map in a $(n-j)$-dimensional polyhedron whose fibres are of diameter less than $2(n-j+1)^{-1/p}$.
\end{proof}
\indent The above proof for an upper bound also gives that $\wdm_\eps (B^{l^p(n)}, l^q) \leq k$ if $\eps \geq 2(k+1)^{1/q-1/p}$, but the inclusion of a $l^q$ ball of proper radius in the $l^p$ ball gives a lower bound that does not meet these numbers. Also note that lemma \ref{kitu} is efficient to evaluate width of tori, as the filling radius of a product is the minimum of the filling radius of each factor.

\section{Evaluation of $\wdm_\eps B^{l^p(n)}_1$} \label{spt2}
We now focus on the computation of $\wdm_\eps X$ for unit ball in finite dimensional $l^p$. Except for a few cases, the complete description is hard to give. We start with a simple example.
\begin{ex} \label{boull1}
Let $B^{l^1(2)}$ be the unit ball of $\rr^2$ for the $l^1$ metric, then
\[
\wdm_\eps B^{l^1(2)}= \left\{
  \begin{array}{llc}
    0  & \textrm{if} &                       \eps \geq 2, \\
    2  & \textrm{if} &                       \eps < 2.
  \end{array} \right.
\]
If $B^{l^1(2)}$ is endowed with the $l^p$ metric, then $\eps < 2^{1/p} \imp \wdm_\eps B^{l^1(2)}= 2$.
\end{ex}
\begin{proof} 
Given any three points whose convex hull contains the origin, two of them have to be on opposite sides, which means their distance is $2^{1/p}$ in the $l^p$ metric. Hence a radial projection is possible for simplices whose vertices form sets of diameter less than $2^{1/p}$. Invoking lemma \ref{wdmineq}, $\frd \del B^{l^1(2)} \geq 2^{-1+1/p}$. Lemma \ref{kitu} concludes. This is specific to dimension $2$ and is coherent with lemma \ref{wdmli}, since, in dimension $2$, $l^\infty$ and $l^1$ are isometric.  
\end{proof}
An interesting lower bound can be obtained thanks to the Borsuk-Ulam theorem; as a reminder, this theorem states that a map from the $n$-dimensional sphere to $\rr^n$ has a fibre containing two opposite points.
\begin{prop} \label{wdmbu}
Let $S= \del B^{l^p(n+1)}_1$ be the unit sphere of a $(n+1)$-dimensional Banach space, then 
\[
\eps <2 \imp \wdm_\eps S > (n-1)/2.
\]
In particular, the same statement holds for $B^{l^p(n+1)}_1$: $\eps <2 \imp \wdm_\eps B^{l^p(n+1)}_1 > (n-1)/2$.
\end{prop}
\begin{proof}
We will show that a map from $S$ to a $k$-dimensional polyhedron, for $k \leq \frac{n-1}{2}$, sends two antipodal points to the same value. Since radial projection is a homeomorphism between $S$ and the Euclidean sphere $S^n = \del B^{l^2(n+1)}_1$ that sends antipodal points to antipodal points, it will be sufficient to show this for $S^n$. Let $f:S^n \to K$ be an $\eps$-embedding, where $K$ is a polyhedron, $\dim K = k \leq (n-1)/2$ and $\eps<2$. Since any polyhedron of dimension $k$ can be embedded in $\rr^{2k+1} $, $f$ extends to a map from $S^n$ to $\rr^n$ that does not associate the same value to opposite points, because $\eps<2$. This contradicts Borsuk-Ulam theorem. The statement on the ball is a consequence of the inclusion of the sphere.
\end{proof}
Hence, $\wdm_\eps B^{l^p(n)}_1$ always jumps from $0$ to at least $\lfloor \frac{n}{2} \rfloor $ if they are equipped with their proper metric.

\paragraph{A first upper bound.} Though this first step is very encouraging, a precise evaluation of $\wdm$ can be convoluted, even for simple spaces. It seems that describing an explicit continuous map with small fibers remains the best way to get upper bounds. Denote by $\id{n}=\{0,\ldots,n\}$.
\begin{lem}\label{proj1}
Let $B$ be a unit ball in a normed $n$-dimensional real vector space. Let $\{p_i\}_{0 \leq i \leq n}$ be points on the sphere $S=\del B$ that are not contained in a closed hemisphere. Suppose that $\forall A \subset \id{n}$ with $\abs{A} \leq n-2$, and $\forall \lambda_j \in \rr_{>0}$, where $j \in \id{n}$, if $\nr{\sum_{i \in A} \lambda_i p_i} \leq 1$, $k \notin A$ and $\nr{\sum_{i \in A} \lambda_i p_i - \lambda_k p_k} \leq 1$, then $\nr{\lambda_k p_k} \leq 1$. A set $p_i$ satisfying this assumption gives
\[
\eps \geq \diam \{p_i\} := \maxx{i \neq j} \nr{p_i-p_j} \imp \wdm_\eps B \leq n-1.
\]
\end{lem}
\begin{proof}
This will be done by projecting the ball on the cone with vertex at the origin over the $n-2$ skeleton of the simplex spanned by the points $p_i$. Note that $n+1$ points satisfying the assumption of this lemma cannot all lie in the same open hemisphere, however we need the stronger hyptothesis that they do not belong to a closed hemisphere. Now let $\Delta_n$ be the $n$-simplex given by the convex hull of $p_0, \ldots, p_n$. We will project the ball on the various convex hulls of $0$ and $n-1$ of the $p_i$. Call $\jo{E}$ the radial projection of elements of the ball (save the origin) to the sphere, and let, for $A \subset \id{n}$, $P_A=\{p_0,\ldots,p_n \} \setminus \{p_i |i \in A \}$. In particular, $P_\vide$ is the set of all the $p_i$. Furthermore, denote by $\ec X$ the convex hull of $X$. Given these notations, $\jo{E}\ec P_{\{i\}}$ is the radial projection of the $(n-1)$-simplex $\ec P_{\{i\}}$ ($\ec P_{\{i\}}$ does not contain $0$ else the points would lie in a closed hemisphere), and $\jo{E} \ec P_{\{i,j\}}$ are parts of the boundary of this projection. Finally, consider, again for  $A\subset \id{n}$, $\Delta'_A = \ec[\jo{E} \ec P_A \cup 0]$. 
\par Let $s_i: \Delta'_{\{i\}} \to \union{ j \neq i }{} \Delta'_{\{i,j\}}$ be the projection along $p_i$. More precisely, we claim that $s_i(p)$ is the unique point of $\Delta'_{\{i,j\}}$ that also belongs to $\Lambda_{p_i}(p)=\{p + \lambda p_i | \lambda \in \rr_{\geq 0} \}$. Existence is a consequence of the fact that the points are not contained in an closed hemisphere, \ie $\exists \mu_i \in \rr_{> 0}$ such that $\sum_{k \in \id{n}} \mu_k p_k =0$. Indeed, $p \in \Delta'_{\{i\}}$, if $p \in \Delta'_{\{i,j\}}$ for some $j$, then there is nothing to show. Suppose that $\forall j\neq i, p \notin \Delta'_{\{i,j\}}$. Then $p = \sum_{k\neq i} \lambda_k p_k$, where $\lambda_k >0$. Write $p_i = -\frac{1}{\mu_i} \sum_{k  \neq i} \mu_k p_k$. It follows that for some $\lambda$, $p + \lambda p_i$ can be written as $\sum_{k \in \id{n} \setminus \{i,j\}} \lambda'_k p_k$ with $0 \leq \lambda'_k \leq \lambda_k$. Uniqueness comes from a transversality observation. $\Delta'_{\{i,j\}}$ is contained in the plane generated by the set $P_{\{i,j\}}$ and $0$ which is of codimension 1. If the line $\Lambda_{p_i}(p)$ was to lie in that plane then the set $P_{\{j\}}$ would lie in the same plane, and $P_\vide$ would be contained in a closed hemisphere. Thus $\Lambda_{p_i}(p)$ is transversal to $\Delta'_{\{i,j\}}$. The figure below illustrates this projection in $\Delta'_{\{0\}}$ for $n=3$.
\vspace{-0.3cm}
\begin{figure}[h]
\centering
\includegraphics{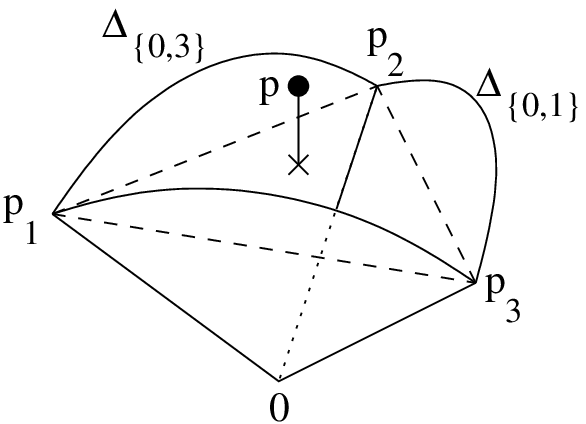}
 \label{toto}
\end{figure}
\vspace{-0.3cm}
\par Our (candidate to be an) $\eps$-embedding $s$ is defined by $s \bv_{\Delta'_{\{i\}}} = s_i$. Since on $\jo{E} \ec P_{\{i\}} \cap \jo{E} \ec P_{\{j\}} \subset \jo{E} \ec P_{\{i,j\}}$, we see that $s\bv_{\Delta'_{\{i,j\}}} = \Id$ and that $\union{i\in \id{n}}{} \Delta'_{\{i \}} = B$, this map is well-defined. It remains to check that the diameter of the fibres is bounded by $\eps$. We claim that the biggest fibre is $s^{-1}(0)= \cup_i \ec \{-p_i,0\}$, whose diameter is that of the set of vertices of the simplex, $\diam \{p_i\}$. To see this, note that for $x \in \Delta'_{\{i,j\}}$, the diameter of $s^{-1}(x)$ attained on its extremal points (by convexity of the norm), that is $x$ and points of the form $x-\lambda_k p_k$ (for $k\in A$, where $A \supset \{i,j\}$ and $x \in \Delta'_A \subset \Delta'_{\{i,j\}}$) whose norm is one. However, since $x = \sum \lambda_i p_i$ for $i \notin A $ and $\lambda_i >0$, $\nr{x-\lambda_k p_k}=1$ implies $\nr{\lambda_k p_k} \leq 1 $, so a simple translation of $s^{-1}(x)$ is actually included in $s^{-1}(0)$.
\end{proof}
This allows us to have a first look at the Euclidean case.
\begin{teo} \label{wdiml2}
Let $B_1^{l^2(n)}$ be the unit ball of $\rr^n$, endowed with the Euclidean metric, and let $b_{n;2} := \sqrt{2(1+\frac{1}{n})}$. Then, for $0<k<n$,
\[
  \begin{array}{rllrll}
      \wdm_\eps B_1^{l^2(n)}&=0&\textrm{if}&        2 \leq & \eps, & \\
k\leq \wdm_\eps B_1^{l^2(n)}&<n&\textrm{if}&  b_{k+1;2} \leq & \eps &< b_{k;2},\\
      \wdm_\eps B_1^{l^2(n)}&=n&\textrm{if}&               & \eps &< b_{n;2}.
  \end{array}
\]
\end{teo}
\begin{proof}
First, when $\eps \geq \diam B_1^{l^2(n)} =2$ this result is a simple consequence of proposition \ref{propbase}.c; when $n=1$ it is sufficient, so suppose from now on that $n\geq 2$. Applying lemma \ref{kitu} to $\del B_1^{l^2(n)} \subset B_1^{l^2(n)}$ yields that $\wdm_\eps B_1^{l^2(n)} =n $ if $\eps < 2 \frd \del B_1^{l^2(n)}$, but $\frd B_1^{l^2(n)} \geq b_{n;2}$ by Jung's theorem (see \cite[\S 2.10.41]{Fed}), as any set whose diameter is less than $<b_{n;2}$ is contained in an open hemisphere (\cite{Kat} shows that $\frd B_1^{l^2(n)}=b_{n;2}$). On the other hand, balls of dimension $k<n$ are all included in $B_1^{l^2(n)}$, which means that $\wdm_\eps B^{l^2(k)}_1 \leq \wdm_\eps B^{l^2(n)}_1$, thanks to \ref{propbase}.d. Hence we have that $\wdm_\eps B_1^{l^2(n)} \geq k$ whenever $b_{k+1;2}\leq \eps < b_{k;2} $. This proves the lower bounds. 
\par The vertices of the standard simplex satisfy the assumption of lemma \ref{proj1}: thanks to the invariance of the norm under rotation we can assume $p_0=(1,0,\ldots,0)$. The other $p_i$ will all have a negative first coordinate, and so will any positive linear combination. Substracting $\lambda p_0$ will be norm increasing. As the diameter of this set is $b_{n;2}$, lemma \ref{proj1} gives the desired upper bound.
\end{proof}
Let us now give an additional upper bound for the $3$-dimensional case:
\begin{prop}\label{projdim3}
If $1\leq p < \infty$, then $\eps \geq 2 (\frac{2}{3})^{1/p} \imp \wdm_\eps B^{l^p(3)}_1 \leq 2  $.
\end{prop}
\begin{proof}
In $\rr^3$ there is a particularly good set of points to define our projections. These are $p_0=3^{\moins\frac{1}{p}}(1,1,1)$, $p_1=3^{\moins\frac{1}{p}}(1,\moins 1,\moins 1)$, $p_2=3^{\moins\frac{1}{p}}(\moins 1,1,\moins 1)$ and  $p_3=3^{\moins\frac{1}{p}}(\moins 1,\moins 1,1)$. Let $x= \lambda_1 p_1$, where $\lambda \in [0,1]$, and suppose $\nr{\lambda_1 p_1 - \lambda_2 p_2}_{l^p} \leq 1$ for $\lambda_2 \in \rr_{\geq 0}$. We have to check that $\lambda_2 \leq 1$. Suppose $\lambda_2 >1$, then $1 \geq \nr{\lambda_1 p_1 - \lambda_2 p_2}_{l^p} = \frac{2}{3}(\lambda_1 + \lambda_2)^p + \frac{1}{3}(\lambda_2-\lambda_1)^p = \lambda_2^p[\frac{2}{3}(1 + t)^p + \frac{1}{3}(1-t)^p]$, where $t= \lambda_1/\lambda_2$. The function of $t$ has minimal value $1$, which gives $\lambda_2 \leq 1 $ as desired. 
\par Suppose now that $x=\lambda_1 p_1 + \lambda_2 p_2$ is of norm less than $1$, where without loss of generality we assume $\lambda_2 \geq \lambda_1$, and $\nr{\lambda_1 p_1 + \lambda_2 p_2 - \lambda_3 p_3}_{l^p} \leq 1$. $\nr{x}_{l^p} \leq 1$ implies that $1 \geq \frac{1}{3}(\lambda_1 + \lambda_2)^p +\frac{2}{3}(\lambda_2 - \lambda_1)^p $ so $(\lambda_2 - \lambda_1)^p \leq 1 - \frac{1}{3}(\lambda_2 + \lambda_1)^p +\frac{1}{3}(\lambda_2 - \lambda_1)^p \leq 1$. If $\lambda_3 > 1$, then 
\[
\begin{array}{rl}
1 &\geq \nr{\lambda_1 p_1 + \lambda_2 p_2 - \lambda_3 p_3}_{l^p} \\
  &  =  \frac{1}{3}(\lambda_3 +\lambda_2 + \lambda_1)^p + \frac{1}{3}(\lambda_3 - (\lambda_2 - \lambda_1))^p +\frac{1}{3}(\lambda_3 +(\lambda_2 - \lambda_1))^p .
\end{array}
\]
However,
\[
\begin{array}{rl}
 \lambda_3^p &\leq  \frac{1}{3}(\lambda_3 +\lambda_2 + \lambda_1)^p + \frac{2}{3}\lambda_3^p \\
             &\leq  \frac{1}{3}(\lambda_3 +\lambda_2 + \lambda_1)^p + \frac{1}{3}(\lambda_3 - (\lambda_2 - \lambda_1))^p +\frac{1}{3}(\lambda_3 +(\lambda_2 - \lambda_1))^p  \\
             &\leq 1.
\end{array}
\]
Using that $f(t) = (1+t)^p + (1-t)^p$ has minimum $2$ for $t \in [0,1]$. These arguments can be repeated for any indices to show that the points $p_i$, where $i =0,1,2$ or $3$, satisfy the assumption of lemma \ref{proj1}. The conclusion follows by showing that $\diam(p_i) = 2 (\frac{2}{3})^{1/p} $
\end{proof}
For certain dimensions, a set of points that allows to build projections with small fibers can be found. Their descriptions require the concept of Hadamard matrices of rank $N$; these are $N\times N$ matrices, that will be denoted $H_N$, whose entries are $\pm 1$ and such that $H_N \cdot H_N^t=N \Id $. It has been shown that they can only exist when $N=2$ or $4|N$, and it is conjectured that this is precisely when they exist. Up to a permutation and a sign, it is possible to write a matrix $H_N$ so that its first column and its first row consist only of $1$s. It is quite easy to see that two rows or columns of such a matrix have exactly $N/2$ identical elements.
\begin{df}\label{hadset}
Let $H_N$ be a Hadamard matrix of rank $N$, and let, for $0 \leq i \leq N$, $h_i$ be the $i^{\textrm{th}}$ row of the matrix without its first entry (which is a $1$). Then the $h_i$ form a Hadamard set in dimension $N-1$.
\end{df}
These $N$ elements, normalised so that $\nr{h_i}_{l^p(N-1)}=1$. When so normalised, their diameter (for the $l^p$ metric) is $2^{1-1/p} (1+\frac{1}{N-1})^p$. Since $\sum h_i=0$, by orthogonality of the columns with the column of $1$ that was removed, we see that they are not contained in an open hemisphere. 
The set of points in the preceding proposition was given by a Hadamard matrix of rank $4$, and when $p=2$ the convex hull of these points is just the standard simplex.
\begin{prop}
Suppose there exists a Hadamard matrix of rank $n+1$, 
then
\[
\eps \geq 1+\frac{1}{n} \imp \wdm_\eps B^{l^1(n)}_1 \leq n-1 .
\]
\end{prop}
\begin{proof}
Let the $h_i$ be as above, and $N=n+1$. Note that for $i\neq j$, $h_i$ and $h_j$ have $\frac{N}{2}$ opposed coordinates, and $\frac{N}{2} -1$ identical ones. Thus $\lambda_i h_i - \lambda_j h_j$ has always a bigger $l^1$ norm than any of its two summands. 
Indeed, the coefficients $c_{j}$ of the vector $\somme{i \in A}{} \lambda_i h_i$ where the contribution of $h_k$ reduces $\abs{c_{j}}$ are in lesser number than those that get increased. Since the $l^1$ norm is linear, the magnitude of the $c_j$ getting smaller is not relevant, only their number.
\par We conclude by applying lemma \ref{proj1}, as $\diam_{l^1}(h_i) = 1 + \frac{1}{N-1}$.
\end{proof}
Note that in dimension higher than $3$ and for $p>2$, Hadamard sets no longer satisfy the assumption of lemma \ref{proj1}.

\paragraph{Further upper bounds for $\wdm_\eps B^{l^p(n)}_1$.}
The projection argument still works for non-Euclidean spheres. It can also be repeated, though unefficiently, to construct maps to lower dimensional polyhedra.
\begin{prop}\label{proj}
For $1<p<\infty$, consider the sphere $B^{l^p(n)}_1$ with its natural metric. Then, for $\frac{n-1}{2} < k < n$, $\exists c_{k,n;p} \in [1,2)$ such that $c_{k,n;p} \geq c_{k+1,n;p}$, and
\[
\wdm_\eps B_1^{l^p(n)} \leq k \textrm{ if }  \eps \geq c_{k,n;p}.
\]
Furthermore $c_{n-1,n;2}=b_{n;2} $ 
\end{prop}
\begin{proof} 
This proposition is also obtained by constructing explicitly maps that reduce dimension (up to $n-j$ for $j<\frac{n+1}{2}$) and whose fibres are small. Unfortunately, nothing indicates this is optimal, and the size of the preimages is hard to determine. We will abbreviate $B:=B_1^{l^p(n)}$.
\par We proceed by induction, and keep the notations introduced in the proof of lemma \ref{proj1}. The $p_i$ that are used here are the vertices of the simplex; they need to be renormalised to be of $l^p$-norm $1$, but note that multiplying them by a constant has actually no effect in this argument. Also note that the sets $\Delta'_A$ are not the same for different $p$, since they are constructed by radial projection to different spheres. The keys to this construction are the maps $s_{j;\{i_1,\ldots,i_j\}}: \Delta'_{\{i_1,\ldots,i_j\}} \to \union{m \notin \{i_1, \ldots, i_j \} }{} \Delta'_{\{i_1, \ldots, i_j,m\} }$ given by projection along the vectors $\somme{l=1}{j} p_{i_l}$. Call $\sigma_1$ the function $s$ from lemma \ref{proj1}, then, for $j>1$, $\sigma_j : B  \to \union{\{i_1,\ldots , i_{j+1} \} \subset \id{n}}{} \Delta'_{\{i_1,\ldots, i_{j+1}\}}$ is obtained by composing, on appropriate domains, $s_{j;\{i_1, \ldots, i_j\}}$ with $\sigma_{j-1}$. Since $s_{j;i_1, \ldots, i_j}$ are equal to the identity when their domain intersect, and their union covers the image of $\sigma_{j-1}$, the map is again well-defined. It remains only to calculate the diameter of the fibres. At $0$ the fibre is
\[ 
\sigma_j^{-1}(0) = \union{\{i_1,\ldots , i_j \} \subset \id{n} }{} \{ -(\lambda_1 + \ldots + \lambda_j)p_{i_1} - (\lambda_1 + \ldots + \lambda_{j-1})p_{i_2} - \ldots - \lambda_1 p_{i_j} | \lambda_i \in \rr_{\geq 0} \}.
\] 
Whereas for a given $x \in \Delta'_A$ in the image (that is $A$ contains at least $j$ elements), $x$ can also be written down as a combination $\sum \lambda_i p_i$, for $i \notin A$ and $\lambda_i \in \rr_{> 0}$. We have
\[ 
\sigma_j^{-1}(x) = \union{\{i_1,\ldots , i_j \} \subset A  }{} \{x -(\lambda_1 + \ldots + \lambda_j)p_{i_1} - (\lambda_1 + \ldots + \lambda_{j-1})p_{i_2} - \ldots - \lambda_1 p_{i_j} | \lambda_i \in \rr_{\geq 0} \}.
\] 
If we set $c_{k,n} = \supp{x \in \sigma_j(B)} \diam \sigma_{n-k}^{-1}(x)$, then when $\eps \geq c_{k,n}$, $\wdm_\eps B^{l^2(n)}_1 \leq k$. It is possible to determine two simple facts about these numbers. First, they are non-increasing $c_{k,n} \geq c_{k+1,n}$, which is obvious as the construction is done by induction, the size of the fiber of maps to lower dimension is bigger than for maps to higher dimension. 
\par Second, they are meaningful: $c_{k,n}<2$. Indeed, when $p \neq 1,\infty$, $c_{k,n}=2$ only if $\sigma^{-1}_{n-k}(x) $ contains opposite points, which is a linear condition.  When $x\neq 0$, by convexity of the distance, the points on which the diameter can be attained are at the boundary of $\sigma_j^{-1}(x)$. Say $Y$ is the set of those point except $x$. The distance from $Y$ to $x$ is at most one, while the diameter of $Y$ is bounded. Indeed, there is a cap of diameter less than $2$ that contains all the $p_i$ but one. The biggest diameter of such caps is also less than $2$ and bounds $\diam Y$.
\par Any point of the fibre at $0$ is a linear combination of the vertices $p_i$, and there is only one linear relation between these, namely $\sum p_i =0$. As long as $j<\frac{n+1}{2}$ (\ie $k > \frac{n-1}{2} $) there are not enough $p_i$ in any two sets that form $\sigma_j^{-1}(0)$ to combine into the required relations, but as soon as $j$ exceeds this bound, opposite points are easily found.
\end{proof}
\par For $B^{l^p(n)}_1$, where $1<p <\infty$, we used the regular simplex to describe our projections, though nothing indicates that this choice is the most appropriate. In fact, many sets of $n+1$ points allow to build projections to a polyhedron, but it is hard to tell which are more effective: on one hand we need this set to have a small diameter (so that the fibre at $0$ is small), while on the other, we need it to be somehow well spread (so as to avoid fibres at $x$ to be too large, as in the assumption of lemma \ref{proj1}). Furthermore, there is in general no reason for $c_{n-1,n;p}$ to coincide with a lower bound, or even to be different from other $c_{k;p}$, thus we cannot always insure that $n-1 \in \wsp (B^{l^p(n)}_1, l^p) $. 

\paragraph{The lowest non zero element of $\wsp$.}
Before we return to the general $l^p$ case, notice that together proposition \ref{wdmbu} and theorem \ref{wdiml2} give a good picture of the function $\wdm_\eps B^{l^2(n)}_1$. It equals $n$ for $\eps < b_{n;2} = c_{n-1,n;2}$, then $n-1$ for $b_{n;2} \leq \eps < b_{n-1;2}$. Afterwards, I could not show a strict inequality for the $c_{k,n;2}$, but even if they are all equal, $\wdm_\eps B^{l^2(n)}_1$ takes at least one value in $(\frac{n}{2}-1, \frac{n}{2}+1) \cap \zz$. Then ,when $\eps \geq 2$, it drops to $0$. 
\par For odd dimensional balls, there is a gap between the value given by proposition \ref{wdmbu} and the lowest dimension obtained by the projections introduced above. Say $B$ is of dimension $2l+1$ and $\eps$ less than but sufficiently close to $2$, then on one hand we know that $\wdm_\eps B \geq l$, while on the other $\wdm_\eps B \leq l+1$. It is thus worthy to ask whether one of these two methods can be improved, perhaps by using extra homological information on the simplices in the proof of proposition \ref{wdmbu} (\eg if its highest degree cohomology is trivial then a $k$-dimensional polyhedron is embeddable in $\rr^{2k}$, see \cite{Fenn}).
\begin{rem}\label{budim3}
Such an improvement is actually available when $n=3$: if the $2$-dimensional sphere maps to a $1$-dimensional polyhedron (\ie a graph), the map lifts to the universal cover, a tree $K$. Hence $K$ is embeddable in $\rr^2$, and, for $1<p<\infty$.
\[
\eps < 2 \imp \wdm_\eps B^{l^p(3)}_1 \geq 2
\]
for otherwise it would contradict Borsuk-Ulam theorem.
\end{rem}
\par Note that estimates obtained in \cite[app 1.E5]{Gfrm} for $\diam_1$, can also yield lower bounds for the diameter of fibres for maps to graphs (\ie $1$-dimensional polyhedra). Applied to spheres, it becomes a special case of proposition \ref{wdmbu} and of the above remark.

\paragraph{Lower bounds for $\wdm_\eps B^{l^p(n)}_1$.}
The remainder of this section is devoted to the improvement of lower bounds, using an evaluation of the filling radius as a product of lemma \ref{wdmineq}, and a short discussion of their sharpness.
\par We shall try to find a lower bound on the diameter of $n+1$ points on the $l^p$ unit sphere that are not in an open hemisphere; recall that points $f_i$ are not in an open hemisphere if $\exists \lambda_i$ such that $\sum \lambda_i f_i =0$. A direct use of Jung's constant (defined as the supremum over all convex $M$ of the radius of the smallest ball that contains $M$ divided by $M$'s diameter) that is cleverly estimated for $l^p$ spaces in \cite{PI} does not yield the result like it did in the Euclidean case. This is due to the fact that there are sets of $n+1$ points on the sphere that are not contained in an open hemisphere, but are contained in a ball (not centered at the origin) of radius less than $1$. The set of points given by
\[ \eqtag \label{contexli}
(1,\ldots ,1), \left(-\frac{2}{n-1},\ldots ,-\frac{2}{n-1},1\right), \ldots , \textrm{ and } \left(1,-\frac{2}{n-1},\ldots ,-\frac{2}{n-1}\right)
\] 
is such an example for $l^\infty$, and deforming it a little can make it work for the $l^p$ case, $p$ finite but close to $\infty$. However, a very minor adaptation of the methods given in \cite{PI} is sufficient.
\par First, we introduce norms for the spaces of sequences (and matrices) taking values in a Banach space $E$. Let $\alpha_i \in \rr_{\geq 0}$ be such that $\somme{i=0}{n} \alpha_i =1$ and denote by $\alpha$ this sequence of $n+1$ real numbers. Let $E_{p,\alpha}$ be the space of sequences made of $n+1$ elements of $E$ and consider the $l^p$ norm weighted by $\alpha$: $\nr{x}_{E_{p,\alpha}} = \big( \sum_i \alpha_i \nr{x_i}_E^p \big)^{1/p}$ where $x=(x_0, \ldots, x_n)$. On the other hand, $E_{p,\alpha^2}$ shall represent the space of matrices whose entries are in $E$, with the norm $\nr{(x_{i,j})}_{E_{p,\alpha^2}} = \big( \sum_{i,j} \alpha_i \alpha_j \nr{x_{i,j}}_E^p \big)^{1/p}$. Now define, for $E,E'$ Banach spaces based on the same vector space and for $1 \leq s,t \leq \infty $, the linear operator $T:E_{s,\alpha} \to E'_{t,\alpha^2}$ by $(x_i) \mapsto (x_i - x_j)$.
\begin{teo} \label{pi1}
Consider a vector space on which two norms are defined, and denote by $E_1$, $E_2$ the Banach space they form. Let $f_i \in E_1^*$, $0 \leq i \leq n$, be such that $\nr{f_i}_{E^*_1}=1$ but that they are not included in an open hemisphere, \ie there exists $\lambda_i \in \rr_{\geq 0} $ such that $\sum \lambda_i f_i =0 $ and $\sum \lambda_i=1 $. Let $\diam_{E^*_2}(f) = \supp{0 \leq i,j \leq n} \nr{f_i-f_j}_{E^*_2}$ be the diameter of this set with respect to the other norm. Then, for $\alpha_i=\lambda_i $,
\[
\diam_{E_2^*}(f) \geq 2 \supp{1\leq s,t \leq \infty} \left(1+ \frac{1}{n} \right)^{1/t'} \supp{E_1} \nr{T}^{-1}_{(E_1)_{s,\alpha} \to (E_2)_{t,\alpha^2}}.
\]
\end{teo}
\begin{proof}
As the $f_i$ are not in an open hemisphere, real numbers $\lambda_i \in \rr_{\geq 0} $ such that $\sum \lambda_i =1$ and $\sum \lambda_i f_i =0$ exist. Furthermore, since $\nr{f_i}_{E_1^*}=1$, there also exist $x_i \in E_1$ such that $f_i(x_i) = 1$ and $\nr{x_i}_{E_1}=1$. The remark on which the estimation relies is, as in \cite{PI},
\[ 
2 = \somme{i,j=0}{n} \lambda_i \lambda_j (f_i-f_j)(x_i-x_j).
\]
Choosing $\alpha_i = \lambda_i$, this equality can be rewritten in the form $2= (Tf)(Tx)$, where $Tx \in (E_2)_{t,\alpha^2}$ and $Tf \in ((E_2)_{t,\alpha^2})^* = (E_2^*)_{t',\alpha^2}$, and thus $2 \leq \nr{Tf}_{(E_2^*)_{t',\alpha^2}} \nr{Tx}_{(E_2)_{t,\alpha^2}}$. Notice that
\[
\somme{i\neq j }{} \alpha_i \alpha_j = \somme{i=0}{n} \alpha_i (1-\alpha_i) \leq 1 - \frac{1}{n+1} = \frac{n}{n+1},
\]
because $\nr{\alpha_i}_{l^1(n+1)} =1 \imp \nr{\alpha_i}_{l^2(n+1)} \geq (n+1)^{-1/2} $. We can isolate the required diameter:
\[
\begin{array}{rl}
\nr{Tf}_{(E_2^*)_{t',\alpha^2}} &= \big( \somme{i=0}{n} \alpha_i \alpha_j \nr{f_i-f_j}_{E_2^*}^{t'} \big)^{1/t'} \\
                                &\leq \diam_{E_2^*}(f) \big( \somme{i\neq j}{} \alpha_i \alpha_j  \big)^{1/t'} \\
                                &\leq  \diam_{E_2^*}(f) \big( \frac{n}{n+1}   \big)^{1/t'} .
\end{array}
\]
On the other hand, $\nr{x_i}_{E_1} =1$, consequently $\nr{x}_{(E_1)_{s,\alpha}} =1 $, so we bound
\[
\nr{Tx}_{(E_2)_{t,\alpha^2}} \leq \nr{T}_{(E_1)_{s,\alpha} \to (E_2)_{t,\alpha^2}}.
\]
The conclusion is found by substitution of the estimates for the norms of $Tf$ and $Tx$.
\end{proof}
We only quote the next result, as there is no alteration needed in that part of the argument of Pichugov and Ivanov. 
\begin{teo} (\cf \cite[thm 2]{PI})
\[
\begin{array}{ll}
\mathrm{if}\; 1\leq p\leq 2      ,&\nr{T}_{(l^p(n))_{\infty,\alpha} \to (l^p(n))_{p,\alpha^2}} \leq 2^{1/p} \left( \frac{n}{n+1} \right)^{1/p - 1/p'},\\
\mathrm{if}\; 2\leq p\leq \infty ,&\nr{T}_{(l^p(n))_{\infty,\alpha} \to (l^p(n))_{p,\alpha^2}} \leq 2^{1/p'}.
\end{array}
\]
\end{teo}
A simple substitution in theorem \ref{pi1}, with $E_1=E_2=l^p(n)$, $s= \infty $ and $t=p$, yields the desired inequalities.
\begin{cor} Let $f_i$, $0\leq i \leq n $, be points on the unit sphere of $l^p(n)$ that are not included in an open hemisphere, then
\[
\begin{array}{ll}
\mathrm{if}\; 1\leq p\leq 2      ,&\diam_{l^p(n)}(f) \geq 2^{1/p'} \left(1+ \frac{1}{n} \right)^{1/p}, \qquad (*)\\
\mathrm{if}\; 2\leq p\leq \infty ,&\diam_{l^p(n)}(f) \geq 2^{1/p} \left(1+ \frac{1}{n} \right)^{1/p'}. \qquad (**)\\
\end{array}
\]
\end{cor}
\begin{rem}
Before we turn to the consequences of this result on $\wdm_\eps$, note that there are examples for which the first inequality is attained. These are the Hadamard sets defined in \ref{hadset}. When normalised to $1$, they are not included in an open hemisphere and of the proper diameter. Hence, when a Hadamard matrix of rank $n+1$ exists, then $(*)$ is optimal. 
Nothing so conclusive can be said for other dimensions, see the argument in example \ref{boull1}. I ignore if there are cases for which $(**)$ is optimal, though it is very easy to construct a family $F_n \in (B^{l^p(n)}_1)^{n+1}$ such that $\diam F_n \to 2^{1/p}$ as $n \to \infty $. In particular for $p=\infty$, the points given in \eqref{contexli} but by substituting $\frac{-1}{n-1}$ instead of the entries with value $\frac{-2}{n-1}$, is a set that is not contained in an open hemisphere and whose diameter is $\frac{n}{n-1}$, which is close to the bound given. Somehow, this case, is also the one where the use of lemma  \ref{wdmineq} results in a bound that is quite far from the right value of $\wdm$, \cf lemma \ref{wdmli}. This might not be so surprising as sets with small diameter on $l^p$ balls seem, when $p>2$, to differ from sets satisfying the assumption of lemma \ref{proj1}.
\end{rem}
\par Still, by lemma \ref{wdmineq} we obtain the following lower bounds on $\wdm$:
\begin{cor} \label{wdimlp}
  Let $b_{k;p}$ be defined by $b_{k;p}= 2^{1/p'} \left(1+ \frac{1}{k}\right)^{1/p}$ when $1\leq p \leq 2$, whereas $b_{k;p} = 2^{1/p} \left(1+ \frac{1}{k}\right)^{1/p'}$ if $2\leq p < \infty$. Then, for $0<k \leq n$,
\[
\eps < b_{k;p} \imp \wdm_\eps B^{l^p(n)}_1 \geq k. 
\]
\end{cor}
\begin{proof}
Let $Y = \del B^{l^p(n)}_1$. Since the convex hull of a set of $n+1$ points on the sphere $Y$ will not contain the origin if the diameter of the set is larger than $b_{n;p}$, lemma \ref{wdmineq} ensures that $\frd Y \geq b_{n;p}/2$. We then use lemma \ref{kitu} for $Y$ to conclude.
\end{proof}
These inequations might not be optimal, proposition \ref{wdmbu} for example is always stronger when $k< \lfloor \frac{n}{2} \rfloor$.
\par In dimension $n$, $B^{l^\infty(n)}_{n^{-1/p}} \subset B^{l^p(n)}_1$ yields that $\eps < 2n^{-1/p} \imp \wdm_\eps B^{l^p(n)}_1 = n $ which improves corollary \ref{wdimlp} as long as 
\[
p \geq \frac{\ln(\frac{2n^2}{n+1})}{\ln (\frac{2n}{n+1})}.
\]
However, when $p=1$, and $H_{n+1}$ is a Hadamard matrix, these estimates are as sharp as we can hope since the lower bound meets the upper bounds.
\begin{cor}
Suppose there is a Hadamard matrix of rank $n+1$. Then, for $0\leq k<n$,
\[
\begin{array}{rllrll}
                          \wdm_\eps B_1^{l^1(n)}&=0&\textrm{if}&                  2 \leq & \eps, & \\
\max(\frac{n-1}{2},k)\leq \wdm_\eps B_1^{l^1(n)}&<n&\textrm{if}&  (1+\frac{1}{k+1}) \leq & \eps &< (1+\frac{1}{k}),\\
                          \wdm_\eps B_1^{l^1(n)}&=n&\textrm{if}&                         & \eps &< (1+\frac{1}{n}).
  \end{array}
\]
\end{cor}
 Furthermore, in dimension $3$, lower bounds of corollary \ref{wdimlp} meet upper bounds of proposition \ref{projdim3} when $1 \leq p \leq 2 $. In particular, thanks to remark \ref{budim3}, this gives a complete description of the $3$-dimensional case for such $p$.
 \begin{cor} \label{wdm3}
Let $p \in [1,2]$, then 
\[
\wdm_\eps B^{l^p(3)}_1 = \left\{
\begin{array}{llrll}
0 & \textrm{if } &2                     \leq &\eps, &     \\
2 & \textrm{if } &2 (\frac{2}{3})^{1/p} \leq &\eps & <  2, \\
3 & \textrm{if } &                           &\eps & <  2 (\frac{2}{3})^{1/p}.
\end{array}
\right.
\]
\end{cor}
When $p>2$, all that can be said is that the value of $\eps$ for which $\wdm_\eps B^{l^p(3)}_1 $ drops from $3$ to $2$ is in the interval $[2 (\frac{2}{3})^{1-1/p} ,2 (\frac{2}{3})^{1/p}]$. 
\par This last corollary is special to the $3$-dimensional case, which happens to be a dimension where there exist a Hadamard set, and where the Borsuk-Ulam argument can be improved to rule out maps to $\frac{n-1}{2}$-dimensional polyhedra. For example, in the $2$-dimensional case, a precise description is not so easy. Indeed, thanks to example \ref{boull1} and using the inclusion of $B^{l^1(2)}_1 \subset B^{l^p(2)}_1$, we know that $\wdm_\eps B^{l^p(2)}_1 = 2$ when $\eps \geq 2^{1/p}$. On the other hand, the inclusion of $B^{l^\infty(2)}_{2^{-1/p}} \subset B^{l^p(2)}_1$ gives $\eps \geq 2^{1/p'} \imp \wdm_\eps B^{l^p(2)}_1 = 2$. Putting these together yields:
\[
\eps \geq \max(2^{1/p},2^{1/p'}) \imp \wdm_\eps B^{l^p(2)}_1 = 2.
\]
These simple estimates in dimension $2$ are better than corollary \ref{wdimlp} as long as $p \leq 3 - \frac{\ln 3}{\ln 2}$ or $p \geq \ln (\frac{8}{3}) / \ln (\frac{4}{3})$. I doubt that any of these estimations actually gives the value of $\eps$ where $\wdm_\eps B^{l^p(2)}_1$ drops from $2$ to $1$.
\par All the results of this section can be summarised to give theorem \ref{t2}. Here are two depictions of the situation. Gray areas correspond to possible values, full lines to known values and dotted line to bounds. 
\vspace{-0.3cm}
\begin{figure}[H]
\centering
\includegraphics{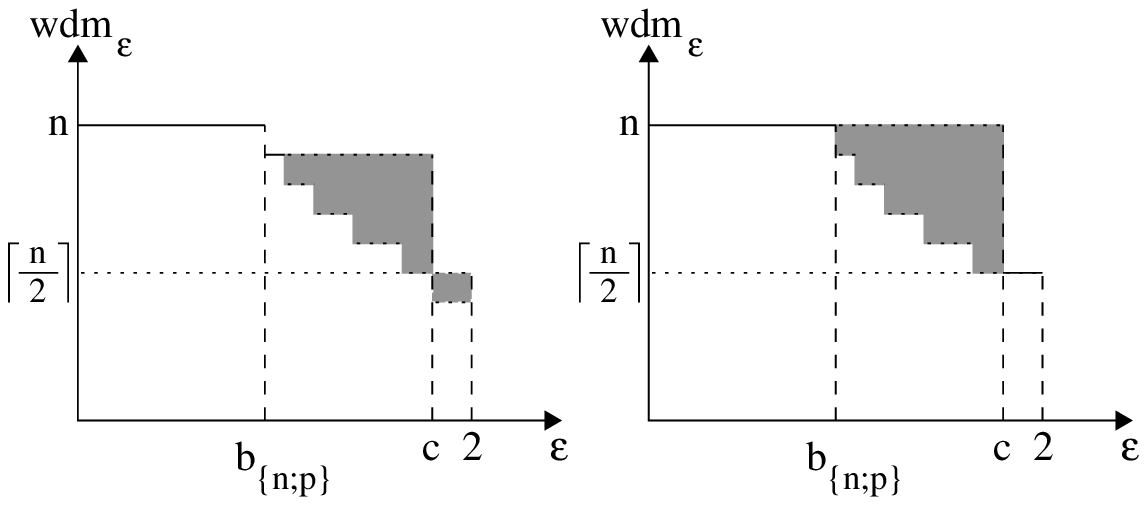}
\end{figure}
\vspace{-0.3cm}
The left-hand plot is for euclidean case, or the case where $p=1$ and there is a Hadamard set, when the dimension is odd and different from $3$: this is when a map to a $n-1$-dimensional polyhedron with small fibers can be constructed, but the bounds from the Borsuk-Ulam argument and projections to lower dimensional polyhedron do not meet. The right-hand one represents the situation in cases where the dimension is even and there is no known projection with small fibers. $c_{\lceil n /2 \rceil, n; p}$ is abbreviated by $c$. The case of dimension $3$ is described in corollary \ref{wdm3}.
\par It is not expected that $\frac{n-1}{2}$ be in $\wsp$ when $n$ is odd, nor is it expected that the lower bounds $b_{k;p}$ be sharp for $B^{l^p(n)}_1$ when $k<n$.


\bibliographystyle{plain}
\bibliography{WidthB}
\end{document}